\documentclass[12pt]{amsart}
\usepackage{amscd,amssymb,amsmath}

\input{prepictex.tex}
\input{pictex.tex}
\input{postpictex.tex}

\def\C{\mathbb C}
\def\N{\mathbb N}
\def\R{\mathbb R}

\def\Z{\mathbb Z}

\def\B{\mathcal B}
\def\H{\mathcal H}
\def\K{\mathcal K}

\def\clsp{\overline{\operatorname{span}}}

\def\id{\operatorname{id}}
\def\ind{\operatorname{ind}}
\def\pr{\operatorname{pr}}

\def\lra{\longrightarrow}
\def\ra{\rightarrow}

\def\irr{\operatorname{Irr}}

\def\1{\rho_{1}^{\theta}}
\def\2{\rho_{2}^{\theta}}
\def\3{\rho_{3}^{\theta}}

\def\n{\rho_{n}^{\theta}}

\def\5{\eta_{1}^{\theta}}
\def\6{\eta_{2}^{\theta}}
\def\7{\eta_{3}^{\theta}}

\newtheorem{theo}{Theorem}[section]
\newtheorem{coro}[theo]{Corollary}
\newtheorem{lemm}[theo]{Lemma}
\newtheorem{prop}[theo]{Proposition}

\newtheorem{rema}[theo]{Remark}

\begin{document}

\title{Noncommutative Balls and Mirror Quantum Spheres}
\author[Hong]{Jeong Hee Hong$^\dag$}
\address{Mathematics, The University of Newcastle, NSW 2308, Australia} 
\curraddr{Applied Mathematics, Korea Maritime University, 
Busan 606--791, South Korea}
\email{hongjh@hhu.ac.kr}

\author[Szyma\'{n}ski]{Wojciech Szyma\'{n}ski$^\ddag$}
\address{Mathematics, The University of Newcastle, NSW 2308, Australia}
\email{Wojciech.Szymanski@newcastle.edu.au}

\thanks{$\dag$ Supported by the Korea Research Foundation
Grant (KRF-2004-041-C00024). \\
\indent $\ddag$ Partially supported by the KBN grant 1 P03A 036 26, 
the European Commission grant MKTD-CT-2004-509794, and the ARC 
Linkage International Fellowship LX0667294.}

\date{2 August, 2007}

\begin{abstract}
Noncommutative analogues of $n$-dimensional balls are defined
by repeated application of the quantum double suspension
to the classical low-dimensional spaces. In the `even-dimensional'
case they correspond to the Twisted Canonical Commutation Relations of
Pusz and Woronowicz. Then quantum spheres are constructed as double manifolds
of noncommutative balls. Both $C^*$-algebras and polynomial algebras
of the objects in question are defined and analyzed, and their relations with
previously known examples are presented. Our construction generalizes that of
Hajac, Matthes and Szyma\'{n}ski for `dimension 2', and leads to a new class of
quantum spheres (already on the $C^*$-algebra level) in all `even-dimensions'.
\end{abstract}

\maketitle

\addtocounter{section}{-1}

\section{Introduction}

Just as classical spheres appear in variety of contexts, their
quantum analogues may be studied from many a different
perspective. One of the most common strategies is to view them as
homogeneous spaces of compact quantum groups \cite{p,frt,vs,hl}.
In addition to quantum symmetry considerations, homological approach
in the spirit of Connes noncommutative geometry has recently become
prominent. Indeed, examples of quantum spheres have been constructed
via Chern character techniques \cite{cl}. We refer the reader to
\cite{d} for an overview of various constructions of quantum
spheres.

Other noncommutative analogues of classical topological methods
have also been used in the study of quantum manifolds, and quantum
spheres in particular. Among them, noncommutative analogues of the
classical suspension were used explicitly or implicitly by several
authors. Quantum double suspension was applied systematically in \cite{hs1,bct},
and noncommutative Heegaard splitting was used in \cite{m,cm,hms3,bhms}.

The main purpose of the present article is to relate quantum spheres to
noncommutative balls, and to examine them from two other natural topological
perspectives. Firstly, we realize quantum spheres as boundaries of
noncommutative balls. Secondly, we construct quantum spheres by
gluing as `double manifolds' of noncommutative balls. Even though the
latter technique goes back to \cite{s}, only recently has it been used
to produce new examples of `two-dimensional' mirror quantum spheres
\cite{hms2}, and we generalize this approach to `higher dimensions'.

In Section 2, working with arbitrary unital $C^*$-algebras and their generators,
we show how to perform the quantum double suspension operation
not only on the $C^*$-algebra level as in \cite{hs1} but also on the
level of a dense $*$-subalgebra (of polynomial functions). In Sections 3 and 4,
we use this procedure to construct noncommutative balls in all `dimensions'
via repeated application of the quantum double suspension to a point (`even
dimensions') and to a closed interval (`odd dimensions'). In Theorems \ref{b2n}
and \ref{bodd}, we present the resulting algebras in terms of convenient
generators and relations. Remarkably, it turns out that in the `even-dimensional'
case our relations are essentially identical with the Twisted Canonical Commutation
Relations of Pusz and Woronowicz \cite{pw}.

The $C^*$-algebra $C(B^{2n}_q)$ of the noncommutative $2n$-ball is generated
by $n$ elements $z_1,\ldots,z_n$. Their commutation relations imply that
$\sum_{i=1}^nz_iz_i^*\leq1$. Thus, it is natural to consider the quotient
of $C(B^{2n}_q)$ by the ideal generated by $1-\sum_{i=1}^nz_iz_i^*$ as the
algebra of functions on the boundary $\partial B^{2n}_q$
of this noncommutative ball. In fact, there is a natural identification
of this boundary with the quantum unitary sphere $S^{2n-1}_\mu$. Similar
considerations apply in the `odd-dimensional' case as well,
with the boundary $\partial B^{2n-1}_q$
identified with the Euclidean quantum sphere $S^{2n-2}_\mu$.

In Sections 5 and 6, we construct the noncommutative double manifold $S^n_{q,\beta}$
of $B^n_q$, by gluing two copies of $B^n_q$ along their common boundary $\partial B^n_q$.
On the $C^*$-algebra level, $C(S^n_{q,\beta})$ is defined by the pull-back
$C(B^n_q)\oplus_\beta C(B^n_q)$ over $C(\partial B^n_q)$. This construction
involves the choice of an automorphism $\beta$ of $C(\partial B^n_q)$,
responsible for the identification of the boundaries of the two noncommutative balls.
Polynomial algebras ${\mathcal O}(S^n_{q,\beta})$ are then defined by
a suitable choice of generators inside $C(S^n_{q,\beta})$.

In the `odd-dimensional' case, it turns out that the isomorphism class of the
$C^*$-algebras $C(S^{2n-1}_{q,\beta})$ does not depend on the choice of $\beta$.
Moreover, these glued quantum spheres $S^{2n-1}_{q,\beta}$ can be naturally
identified with the unitary quantum spheres (Proposition \ref{s2n1}). The
situation is quite different in the `even-dimensional' case. Indeed, we find
an automorphism $\beta$ of $C(\partial B^{2n}_q)$
such that the $C^*$-algebra $C(S^{2n}_{q,\beta})$ is not even stably isomorphic
to $C(S^{2n}_{q,\id})$ (Corollary \ref{s2nbetaid}). This happens in spite of
the fact that these two $C^*$-algebras (of type $I$) have homeomorphic primitive
ideal spaces and isomorphic (classical) $K$-groups (Theorem \ref{s2nbetaK}).
For such a $\beta$, we call $S^{2n}_{q,\beta}$ {\em mirror quantum sphere}.
Our construction generalizes that of \cite{hms2} carried for `dimension 2'. While
$S^{2n}_{q,\id}$ may be naturally identified with the Euclidean quantum spheres
(Proposition \ref{s2nid}), the mirror quantum spheres $S^{2n}_{q,\beta}$ are new
(already on the $C^*$-algebra level).

Finally, in Section 7, irreducible representations of the $C^*$-algebras of
noncommutative balls $B^n_q$ and of the mirror quantum spheres $S^{2n}_{q,\beta}$
are presented.

\vspace{2mm}\noindent
{\bf Acknowledgements.} The second named author would like to thank Piotr
Hajac and the entire team of the program in Noncommutative Geometry and
Quantum Groups for their warm hospitality during his stay in Warsaw in
March--May 2006.


\section{The double of a noncommutative space}

Let $X$ be a compact manifold with non-empty boundary $\partial X$.
Given a homeomorphism $f:\partial X\ra\partial X$ of the boundary, the classical
topological gluing construction yields a double $X\cup_f X$ of $X$.
To translate this picture into the language of $C^*$-algebras, let
$C(X)$ be the commutative $C^*$-algebra of continuous complex-valued functions on
$X$, and let $C_{\partial X}(X)$ denote the continuous complex-valued functions on
$X$ vanishing on $\partial X$. Then $C_{\partial X}(X)$ is an essential
ideal of $C(X)$. If $\pi:C(X)\ra C(\partial X)$ is the
surjection given by restriction then we have an exact sequence
of commutative $C^*$-algebras
\begin{equation}\label{commutativeboundary}
0 \lra C_{\partial X}(X) \lra C(X) \overset{\pi}{\lra} C(\partial X) \lra 0.
\end{equation}

The $C^*$-algebra $C(X\cup_f X)$ is isomorphic to the pull-back of
$C(\partial X)$ along two surjections $\pi:C(X)\ra C(\partial X)$ and
$f_*\circ\pi:C(X)\ra C(\partial X)$, where $f_*:C(\partial X)\ra
C(\partial X)$ is the automorphism dual to $f$.

\begin{rema}\label{algebraic}
\rm An imbedding of $X$ into a Euclidean space ($\R^n$ or $\C^n$) gives
rise to a dense $*$-subalgebra ${\mathcal O}(X)$ (the polynomial algebra) of
$C(X)$, generated by the restrictions of the coordinate functions to $X$.
Clearly, we have $\pi({\mathcal O}(X))={\mathcal O}(\partial X)$. Furthermore,
$\partial X$ is the intersection of $X$ with an affine variety if
and only if ${\mathcal O}(X)\cap C_{\partial X}(X)$ is dense in $C_{\partial X}(X)$
in the sup norm topology.
\end{rema}

In the present article, we are concerned with noncommutative
analogues of the aforementioned classical setting.
Let $A$ be a unital $C^*$-algebra (not necessarily commutative),
$J$ be an essential proper ideal of $A$, and $\pi:A\rightarrow B=A/J$
be the natural surjection. Thus, we have an essential extension
\begin{equation}\label{boundaryextension}
0 \lra J \lra A \overset{\pi}{\lra} B \lra 0.
\end{equation}
Suppose that $\beta$ is an automorphism of $B$. Then we define
$A{\oplus}_\beta A$, the double of $A$, as
\begin{equation}\label{double}
A\oplus_\beta A=\{(x,y)\in A\oplus A:\pi(x)=(\beta\circ\pi)(y)\}.
\end{equation}
That is, $A\oplus_\beta A$ is the $C^*$-algebra defined by the
pull-back diagram

\vspace{-2mm}
\[
\begin{CD}
A\oplus_\beta A @> {\pr_2} >> A \\
@V {\pr_1} VV  @VV {\beta\circ\pi} V \\
A @>  {\pi} >> B
\end{CD}
\]

\vspace{1mm}
The following proposition follows at once from our definitions.
\begin{prop}\label{doubletype}
Let $\beta$ and $\gamma$ be such automorphisms of $B$ that
$\gamma^{-1}\beta$ admits a lift to an automorphism $\alpha$ of $A$,
that is $\pi\alpha=(\gamma^{-1}\beta)\pi$. Then the map
$A\oplus_\beta A\ra A\oplus_\gamma A$ given by $(x,y)\mapsto
(x,\alpha(y))$ is an isomorphism between these two $C^*$-algebras.
\end{prop}

Combining (\ref{boundaryextension}) with (\ref{double}) we get the
following exact sequence for the double $A\oplus_\beta A$:
\begin{equation}\label{doubleexact}
0 \lra J\oplus J \lra A\oplus_\beta A \lra B \lra 0,
\end{equation}
which may be used to determine structural properties and invariants
of the double.

\begin{rema}\label{algebraicNC}
\rm Throughout this article, an important role is played not only
by $C^*$-algebras but also by their polynomial subalgebras. Thus,
if $A$ is a $C^*$-algebra, as above, we will also consider a dense
$*$-subalgebra ${\mathcal O}$ of $A$, the polynomial algebra
of $A$. Then $\pi({\mathcal O})$ will play the role of
a polynomial algebra of $B$. In our examples, the following two
additional conditions will be satisfied; firstly, ${\mathcal O}\cap J$ will
be dense in $J$ (cf. Remark \ref{algebraic}), secondly,
the automorphisms $\beta$ of $B$ will be algebraic in the sense that
$\beta(\pi({\mathcal O}))=\pi({\mathcal O})$.
\end{rema}


\section{The quantum double suspension}

Let $\{\xi_j:j=0,1,\ldots\}$ be the usual orthonormal basis of
$\ell^2(\N)$. We denote by $\{E_{ij}:i,j=0,1,\ldots\}$ the system of
rank one matrix units such that $E_{ij}(\xi_k)=\delta_{jk}\xi_i$.
Then $\clsp\{E_{ij}:i,j=0,1,\ldots\}$ coincides with the algebra $\K$
of compact operators. We denote by $V$ the unilateral shift
$V(\xi_j)=\xi_{j+1}$. We also denote by $z=\id$ the standard unitary
generator of $C(S^1)$ with $S^1$ the unit circle in $\C$.

Let $A$ be a unital $C^*$-algebra.
In \cite{hs1}, the quantum double suspension $\Sigma^2A$ of $A$ was defined
as the unital $C^*$-algebra for which there exists an essential extension
\begin{equation}\label{doublesuspension}
0 \lra A\otimes\K \lra \Sigma^2A \lra C(S^1) \lra 0
\end{equation}
whose Busby invariant $C(S^1)\ra M(A\otimes\K)/A\otimes\K$ sends
$z$ to the image of $1\otimes V$ under the natural surjection
$M(A\otimes\K)\ra M(A\otimes\K)/A\otimes\K$. The following proposition
follows easily from the definition of the quantum double suspension and
its proof is omitted.

\begin{prop}\label{qdoublesusp}
The quantum double suspension $\Sigma^2A$ of $A$ is isomorphic with the
$C^*$-subalgebra of $A\otimes\B(\ell^2(\N))$ generated by
$A\otimes E_{00}$ and $1\otimes V$.
\end{prop}

In what follows, we identify $\Sigma^2A$ with the $C^*$-algebra of
Proposition \ref{qdoublesusp}. Using this proposition one easily
derives the following useful universal property of the quantum double
suspension.

\begin{prop}\label{universalqds}
The quantum double suspension $\Sigma^2A$ of $A$ has the following
universal property. If $B$ is a unital $C^*$-algebra, $T$ is an
isometry in $B$ {\em (}that is, $T^*T=1${\em )}, and $\tilde{\psi}:A\ra B$ is a
$*$-homomorphism such that $\tilde{\psi}(1)=1-TT^*$, then there exists
a unique $*$-homomorphism $\psi:\Sigma^2A\ra B$ such that
$\psi(1\otimes V)=T$ and $\psi(a\otimes E_{00})=\tilde{\psi}(a)$
for all $a\in A$.
\end{prop}

We want to work with the quantum double suspension in the
following setting. A unital $C^*$ algebra $A$ plays the role
of a deformed function algebra on a compact manifold. The deformation
parameter is denoted $q$ and belongs to $(0,1)$, with $q=1$
being the classical case. Furthermore, $x_1,\ldots,x_n$ are distinguished
elements of $A$ generating a dense $*$-subalgebra (the polynomial algebra).
There is no canonical way to define a polynomial algebra inside
$\Sigma^2A$ and one has to make a choice. Using the notation introduced
earlier in this section, we select the following generators
$X_1,\ldots,X_n,X_{n+1}$ of the polynomial algebra of $\Sigma^2A$:
\begin{eqnarray}
X_j & = & x_j\otimes\sum_{k=0}^\infty q^{k/2}E_{kk} \label{s2generators1} \\
 & = & \sum_{k=0}^\infty q^{k/2}(1\otimes V)^k(x_j\otimes E_{00})
 (1\otimes V^*)^k,  \;\;\;\;\; \text{for}\; j=1,\ldots,n, \nonumber \\
X_{n+1} & = & 1\otimes\sum_{k=0}^\infty\sqrt{1-q^{k+1}}VE_{kk} \label{s2generators2} \\
 & = & \sum_{m=0}^\infty\left(\sqrt{1-q^{m+1}}-\sqrt{1-q^m}\right)
 (1\otimes V)^{m+1}(1\otimes V^*)^m. \nonumber
\end{eqnarray}
Clearly, $X_1,\ldots,X_{n+1}$ all belong to $\Sigma^2A$ and generate
its dense $*$-subalgebra. Furthermore, the above choice of generators
for $\Sigma^2A$ yields the following identities.
\begin{eqnarray}
X_jX_{n+1} & = & q^{1/2}X_{n+1}X_j, \;\;\; j=1,\ldots,n, \label{s2identities1} \\
X_jX_{n+1}^* & = & q^{-1/2}X_{n+1}^*X_j, \;\;\; j=1,\ldots,n, \label{s2identities2} \\
X_{n+1}^*X_{n+1}-qX_{n+1}X_{n+1}^* & = & 1-q. \label{s2identities3}
\end{eqnarray}

An easy calculation yields the following identity, relating
$\{X_j:j=1,\ldots,n+1\}$ to $\{x_j:j=1,\ldots,n\}$, which will be very
useful later in this paper.
\begin{equation}\label{oddradius}
1-\sum_{j=1}^{n+1}X_jX_j^*=\left(1-\sum_{j=1}^n x_jx_j^*\right)
\otimes\sum_{k=0}^\infty q^kE_{kk}.
\end{equation}

\begin{rema}\label{others2generators}
\rm Other natural choices of generators of $\Sigma^2A$ are also possible.
For example, one might set $X_j=x_j\otimes\sum_{k=0}^\infty q^{k}E_{kk}$,
$j=1,\ldots,n$ in (\ref{s2generators1}), with $X_{n+1}$ as in (\ref{s2generators2}).
Then in formulae (\ref{s2identities1})--(\ref{s2identities2})
we get $q$ and $q^{-1}$ instead of $q^{1/2}$ and $q^{-1/2}$, respectively.
\end{rema}


\section{The `even-dimensional' noncommutative balls}

\subsection{The algebra}

In \cite{hs1}, the $C^*$-algebras $C(B^{2n}_q)$ of continuous functions
on `even-dimensional' noncommutative balls were
defined inductively for all $n$ as follows:
$$  C(B^{0}_q)=\C, \;\;\;\;\; C(B^{2(n+1)}_q)=\Sigma^2C(B^{2n}_q). $$
In this definition the parameter $q$ is not explicitly involved,
but it will become visible with a
suitable choice of a dense $*$-subalgebra (polynomial algebra), below.

To begin with, we describe basic properties of these $C^*$-algebras
(cf. \cite{hs1,ml}). As shown in \cite{hs1},
$C(B^{2n}_q)$ is isomorphic with the $C^*$-algebra
$C^*(M_n)$ of a directed graph we call $M_n$. This graph consists of
$n+1$ vertices $\{v_1,\ldots,v_{n+1}\}$, and edges $\{e_{i,j}:i=1,\ldots,n,\;
j=i,\ldots,n+1\}$ such that the source $s(e_{i,j})$ of $e_{i,j}$ is
$v_i$ and its range $r(e_{i,j})$ is $v_j$. In particular, $M_n$ contains
a unique sink $v_{n+1}$ (vertex emitting no edges), and from every
other vertex there is a path to this one. For example, graph $M_3$
looks as follows.

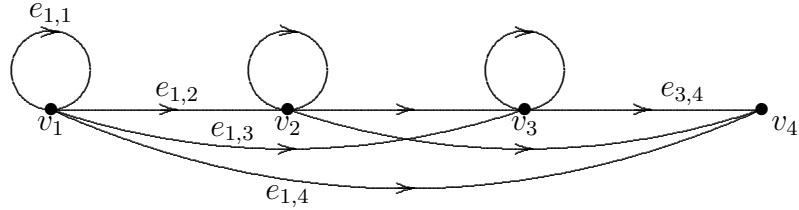
\vspace{1mm}

\[ \beginpicture
\setcoordinatesystem units <2.1cm,2.1cm>
\setplotarea x from -4 to 2, y from -0.7 to 0.6
\put {$\bullet$} at -3 0
\put {$\bullet$} at -1.5 0
\put {$\bullet$} at 0 0
\put {$\bullet$} at 1.5 0
\setlinear
\plot -3 0 -1.5 0 /
\plot -1.5 0 0 0 /
\plot 0 0 1.5 0 /
\arrow <0.235cm> [0.2,0.6] from -2.3 0 to -2.2 0
\arrow <0.235cm> [0.2,0.6] from -0.8 0 to -0.7 0
\arrow <0.235cm> [0.2,0.6] from 0.7 0 to 0.8 0
\setquadratic
\plot -3 0 -0.75 -0.5 1.5 0 /
\plot -3 0 -1.5 -0.25 0 0 /
\plot -1.5 0  0 -0.25 1.5 0 /
\arrow <0.235cm> [0.2,0.6] from -1.5 -0.249 to -1.45 -0.24852
\arrow <0.235cm> [0.2,0.6] from 0 -0.249 to 0.05 -0.24852
\arrow <0.235cm> [0.2,0.6] from -0.8 -0.499 to -0.7 -0.5
\circulararc 360 degrees from -3 0 center at -3 0.25
\circulararc 360 degrees from -1.5 0 center at -1.5 0.25
\circulararc 360 degrees from 0 0 center at 0 0.25
\arrow <0.235cm> [0.2,0.6] from -3 0.492 to -2.95 0.492
\arrow <0.235cm> [0.2,0.6] from -1.49 0.495 to -1.44 0.49
\arrow <0.235cm> [0.2,0.6] from 0 0.49 to 0.05 0.49
\put {$e_{1,1}$} [v] at -3 0.6
\put {$e_{1,2}$} [v] at -2.2 0.1
\put {$e_{1,3}$} [v] at -1.85 -0.15
\put {$e_{1,4}$} [v] at -1.5 -0.55
\put {$e_{3,4}$} [v] at 1 0.1
\put {$v_1$} [v] at -3 -0.1
\put {$v_2$} [v] at -1.5 -0.1
\put {$v_3$} [v] at 0 -0.1
\put {$v_4$} [v] at 1.65 -0.1
\endpicture \]

\vspace{1mm}

Then $C^*(M_n)$ is, by definition, the universal $C^*$-algebra generated by
mutually orthogonal projections $\{P_i:i=1,\ldots,n+1\}$ (corresponding to
the vertices of the graph) and partial isometries $\{S_{i,j}:i=1,\ldots,n,
\;j=i,\ldots,n+1\}$ (corresponding to the edges), subject to the
relations: $S_{i,j}^*S_{i,j}=P_j$ and $P_i=\sum_{j=i}^{n+1}S_{i,j}S_{i,j}^*$.
Now the ideal structure of $C(B^{2n}_q)$ follows from the general
theory of graph algebras \cite{bhrs,hs2}. In particular, its primitive ideal
space consists of $n$ circles and one point.
Similarly, the results of \cite{rsz} yield the
$K$-theory of $C(B^{2n}_q)\cong C^*(M_n)$ as $K_1(C(B^{2n}_q))=0$, and
$K_0(C(B^{2n}_q))\cong\Z$ generated by the class of identity $[1]_0$.

\vspace{2mm}
The following theorem provides a convenient set of generators for
$C(B^{2n}_q)$, which may be regarded as $q$-deformed complex
coordinate functions on a unit ball in $\C^n$.

\begin{theo}\label{b2n}
$C(B^{2n}_q)$ is isomorphic with the $C^*$-algebra $C^*(z_1,\ldots,z_n)$,
universal {\em(}with respect to representations in bounded operators{\em)}
for the relations
\begin{eqnarray}
z_iz_j & = & q^{1/2}z_jz_i  \;\;\;\;\; \mbox{for} \;\; i<j,
\label{b2npresentation1} \\
z_iz_j^* & = & q^{-1/2}z_j^*z_i  \;\;\;\;\; \mbox{for} \;\; i<j,
\label{b2npresentation2} \\
z_i^*z_i-qz_iz_i^* & = & (1-q)\left(1-\sum_{j=i+1}^n z_jz_j^*\right)
\;\;\;\;\; \mbox{for} \;\; i=1,\ldots,n. \label{b2npresentation3}
\end{eqnarray}
\end{theo}
\begin{proof}
At first we observe that the universal norm for relations
(\ref{b2npresentation1})--(\ref{b2npresentation3}) is finite. Indeed,
(\ref{b2npresentation3}) with $i=n$ yields $z_n^*z_n=qz_nz_n^*+1-q$.
Since $||z_n^*z_n||=||z_nz_n^*||=||z_n||^2$ and $q\neq1$, this
gives $||z_n||=1$. Then proceeding by reverse induction on
$i$ we deduce from (\ref{b2npresentation3}) that the norms of all
$z_1,\ldots,z_n$ are universally bounded.

To prove the main part of the theorem, we proceed by induction on $n$.
If $n=1$ then the claim is that $C(B^2_q)$ is isomorphic with the
$C^*$-algebra $C^*(z_1)$, universal for the relation $z_1^*z_1-qz_1z_1^*=1-q$.
The latter is nothing but the Klimek-Lesniewski algebra of the quantum disc,
isomorphic with the Toeplitz algebra \cite{kl}. On the other hand,
$C(B^2_q)\cong C^*(M_1)$ (the quantum double suspension of the complex
numbers) is the graph algebra known to be isomorphic to the Toeplitz algebra.
This establishes the base for induction.

For the inductive step, suppose $C(B^{2n}_q)\cong C^*(z_1,\ldots,z_n)$. Let
$C^*(Z_1,\ldots,Z_{n+1})$ be the universal $C^*$-algebra for
relations (\ref{b2npresentation1})--(\ref{b2npresentation3}), with $n$
replaced by $n+1$. We must show that there exist $*$-homomorphisms
$\phi:C^*(Z_1,\ldots,Z_{n+1})\ra \Sigma^2C^*(z_1,\ldots,z_n)$ and
$\psi:\Sigma^2C^*(z_1,\ldots,z_n)\ra C^*(Z_1,\ldots,Z_{n+1})$ such
that $\psi\circ\phi=\id$ and $\phi\circ\psi=\id$.

To construct $\phi$, we map $Z_1,\ldots,Z_{n+1}$ to the $n+1$ generators
of $\Sigma^2C^*(z_1,\ldots,z_n)$ given by formulae
(\ref{s2generators1})--(\ref{s2generators2}). That is,
\begin{eqnarray}
\phi(Z_i) & = & z_i\otimes\sum_{k=0}^\infty q^{k/2}E_{kk}
\;\;\;\;\; \mbox{for} \;\; i=1,\ldots,n, \label{b2nphi1} \\
\phi(Z_{n+1}) & = & 1\otimes\sum_{k=0}^\infty\sqrt{1-q^{k+1}}VE_{kk}.
\label{b2nphi2}
\end{eqnarray}
It is not difficult to verify that the elements $\phi(Z_1)$,\ldots,
$\phi(Z_{n+1})$ satisfy (\ref{b2npresentation1})--(\ref{b2npresentation3}).
Thus, $\phi$ extends to a $*$-homomorphism from $C^*(Z_1,\ldots,Z_{n+1})$
to $\Sigma^2C^*(z_1,\ldots,z_n)$.

To define $\psi$, we first observe that $Z_{n+1}^*Z_{n+1}=
qZ_{n+1}Z_{n+1}^*+1-q\geq 1-q$ holds by
(\ref{b2npresentation3}). Thus $Z_{n+1}^*Z_{n+1}$ is
invertible and hence $Z_{n+1}$ admits a polar decomposition $Z_{n+1}=
T|Z_{n+1}|$ in $C^*(Z_1,\ldots,Z_{n+1})$, with $T$ an isometry.
We have $(1-TT^*)Z_{n+1}=0$. Since $Z_{n+1}Z_{n+1}^*$ commutes with
$Z_1,\ldots,Z_n$ (by (\ref{b2npresentation1}) and (\ref{b2npresentation2})),
so does $1-TT^*$. Then one can verify that the elements $Z_i(1-TT^*)$,
$i=1,\ldots,n$, satisfy relations (\ref{b2npresentation1})--(\ref{b2npresentation3}).
Consequently, there exists a $*$-homomorphism $\tilde{\psi}:C^*(z_1,\ldots,z_n)\ra
C^*(Z_1,\ldots,Z_{n+1})$ such that $\tilde{\psi}(z_i)=Z_i(1-TT^*)$ for
$i=1,\ldots,n$. Since $\tilde{\psi}(1)=1-TT^*$, Proposition \ref{universalqds}
implies that there is a $*$-homomorphism $\psi:\Sigma^2C^*(z_1,\ldots,z_n)\ra
C^*(Z_1,\ldots,Z_{n+1})$ such that $\psi(1\otimes V)=T$ and
$\psi(a\otimes E_{00})=\tilde{\psi}(a)$ for all $a\in C^*(z_1,\ldots,z_n)$.

We now verify that $\phi\circ\psi=\id$. Indeed, we have
$$
(\phi\circ\psi)(1\otimes V)=\phi(T)=
\phi(Z_{n+1})\phi(Z_{n+1}^*Z_{n+1})^{-1/2}=1\otimes V,
$$
and for $i=1,\ldots,n$ we have
$$
(\phi\circ\psi)(z_i\otimes E_{00})=\phi(Z_i)\phi(1-TT^*)=
\phi(Z_i)(1\otimes E_{00})=z_i\otimes E_{00}.
$$
These two identities imply that $\phi\circ\psi=\id$.
To prove $\psi\circ\phi=\id$ we need the following technical lemma.
\begin{lemm}\label{sumidentities}
Let $Z_{n+1}=T|Z_{n+1}|$ be the polar decomposition in
$C^*(Z_1,\ldots,Z_{n+1})$. Then the following identities hold:
\begin{eqnarray}
Z_i & = & \sum_{k=0}^\infty q^{k/2}T^kZ_i(1-TT^*)(T^*)^k
\;\;\;\;\; \mbox{for} \;\; i=1,\ldots,n,
\label{sumidentities1} \\
Z_{n+1} & = & T+\sum_{k=0}^\infty(\sqrt{1-q^{k+1}}-1)T^{k+1}
(1-TT^*)(T^*)^k.
\label{sumidentities2}
\end{eqnarray}
\end{lemm}
\begin{proof}
Let $\rho$ be a representation of $C^*(Z_1,\ldots,Z_{n+1})$ on
a Hilbert space $\H$. Then, using the Wold decomposition of $\rho(T)$,
$\H$ is a direct sum of two subspaces $\H_0$ and $\H_1$ such that
$\rho(T)|_{\H_0}$ is unitary and $\H_1$ is a direct sum of
the subspaces $\rho(T^k(1-TT^*)(T^*)^k)\H$ for $k=0,1,\ldots$. Thus
$\rho(Z_{n+1})|_{\H_0}$ is invertible and it follows from
(\ref{b2npresentation3}) that it is unitary, and hence
$\rho(Z_{n+1})|_{\H_0}=\rho(T)|_{\H_0}$. Then, using
(\ref{b2npresentation3}) again, one shows by reverse induction on
$j$ that all $\rho(Z_j)|_{\H_0}$ are zero for $j=1,\ldots,n$.
Thus it suffices to verify (\ref{sumidentities1}) and (\ref{sumidentities2})
on $\H_1$, and this follows easily from the identities
(\ref{b2npresentation1})--(\ref{b2npresentation3}).
\end{proof}
Back to the proof of $\psi\circ\phi=\id$. By virtue of
(\ref{sumidentities2}), we have
$$
(\psi\circ\phi)(Z_{n+1})=\psi\left(1\otimes\sum_{k=0}^\infty
\sqrt{1-q^{k+1}}VE_{kk}\right)
$$
$$
=\psi(1\otimes V)+\psi\left(\sum_{k=0}^\infty(\sqrt{1-q^{k+1}}-1)
(1\otimes V^{k+1})(1\otimes(1-VV^*))(1\otimes V^*)^k\right)
$$
$$
=T+\sum_{k=0}^\infty(\sqrt{1-q^{k+1}}-1)T^{k+1}(1-TT^*)(T^*)^k=Z_{n+1}.
$$
Similarly, (\ref{sumidentities1}) implies that for $i=1,\ldots,n$
we have
$$
(\psi\circ\phi)(Z_i)=\psi\left(z_i\otimes\sum_{k=0}^\infty q^{k/2}E_{kk}\right)
=\psi\left(\sum_{k=0}^\infty q^{k/2}(1\otimes V^k)(z_i\otimes E_{00})
(1\otimes V^*)^k\right)
$$
$$
=\sum_{k=0}^\infty q^{k/2}T^kZ_i(1-TT^*)(T^*)^k=Z_i.
$$
These identities entail $\psi\circ\phi=\id$, as required.
\end{proof}
In the remainder of this paper, we suppress the isomorphism $\phi$ of
Theorem \ref{b2n} and simply identify $C(B^{2n}_q)=C^*(z_1,\ldots,z_n)$.
We also note that in the course of proof of this theorem we showed
that generators $z_1,\ldots,z_n$ of $C(B^{2n}_q)$ and
$Z_1,\ldots,Z_{n+1}$ of $C(B^{2(n+1)}_q)$ are related to one another
in accordance with formulae (\ref{s2generators1})--(\ref{s2generators2}),
and hence they satisfy identity (\ref{oddradius}).

We define the polynomial algebra ${\mathcal O}(B^{2n}_q)$ as the $*$-subalgebra of
$C(B^{2n}_q)$ generated by $z_1,\ldots,z_n$. Remarkably, it turns out that
this is exactly the $*$-algebra of twisted canonical commutation relations
(TCCR) of Pusz and Woronowicz \cite{pw}. Indeed, the $*$-algebra of TCCR
is generated by elements $a_1,\ldots,a_n$ satisfying the relations
\begin{eqnarray}
a_ja_i & = & \mu a_ia_j \;\;\;\;\; \mbox{for} \;\; i<j, \label{TCCR1} \\
a_ja_i^* & = & \mu a_i^*a_j \;\;\;\;\; \mbox{for} \;\; i\neq j, \label{TCCR2} \\
a_ia_i^* & = & 1+\mu^2a_i^*a_i-(1-\mu^2)\sum_{j=i+1}^na_j^*a_j
\;\;\;\;\; \mbox{for} \;\; i=1,\ldots,n. \label{TCCR3}
\end{eqnarray}
Setting $q=\mu^2$ one obtains an identification of this algebra
with our ${\mathcal O}(B^{2n}_q)$ through the simple transformation:
\begin{equation}\label{TCCR=b2n}
z_i=\sqrt{1-q}a_i^*, \;\;\;\;\; \mbox{for} \;\; i=1,\ldots,n.
\end{equation}
Consequently, $C(B^{2n}_q)$ is isomorphic with the enveloping $C^*$-algebra
of TCCR. This fact could also be derived from our realization of
$C(B^{2n}_q)$ as the graph algebra $C^*(M_n)$, combined with the
stability results for TCCR obtained in \cite{ps} and \cite{jps}.

\subsection{The boundary}

We first observe that the generators of $C(B^{2n}_q)$ satisfy
\begin{equation}\label{evenradiusinequality}
\sum_{j=1}^n z_jz_j^*\leqq 1.
\end{equation}
This inequality is established by induction on $n$, as follows.
In $C(B^2_q)$, we have $||z_1||=1$ and hence $z_1z_1^*\leqq 1$.
Assuming $\sum_{j=1}^n z_jz_j^*\leqq 1$ in $C(B^{2n}_q)$, we derive
$\sum_{j=1}^{n+1}Z_jZ_j^*\leqq 1$ in $C(B^{2(n+1)}_q)=\Sigma^2C(B^{2n}_q)$
from identity (\ref{oddradius}).

We define $J_{2n}$ as the closed two-sided ideal of $C(B^{2n}_q)$ generated by
$1-\sum_{j=1}^n z_jz_j^*$.

\begin{lemm}\label{J2n}
$J_{2n}$ is an essential ideal of $C(B^{2n}_q)$. Furthermore, $J_{2n}$
is isomorphic with the compacts $\K$ and satisfies $J_{2(n+1)}=J_{2n}\otimes\K$.
\end{lemm}
\begin{proof}
To prove the lemma, we use the identification of $C(B^{2n}_q)$ with
the graph algebra $C^*(M_n)$. The general theory of graph algebras
tells us that the closed two-sided ideal $\tilde{J}_{2n}$ of $C^*(M_n)$
generated by projection $P_{n+1}$ is essential and isomorphic to the compacts
$\K$ \cite{dhs}. Moreover, it is the only ideal of $C^*(M_n)$ with these
two properties. An easy inductive argument shows that $\tilde{J}_{2n}
\otimes\K$ is an essential ideal of $C(B^{2(n+1)}_q)$ isomorphic with the
compacts. Hence we have $\tilde{J}_{2(n+1)}=\tilde{J}_{2n}\otimes\K$.

Therefore, it is enough to verify that $\tilde{J}_{2n}$ coincides with $J_{2n}$, that is
$1-\sum_{j=1}^n z_jz_j^*$ generates $\tilde{J}_{2n}$. For this it suffices to
check that $1-\sum_{j=1}^n z_jz_j^*$ belongs to $\tilde{J}_{2n}$,
since $\tilde{J}_{2n}$ is isomorphic with $\K$. We
proceed by induction on $n$. If $n=1$ then $C(B^{2}_q)$ is isomorphic
with the Toeplitz algebra and $1-z_1z_1^*$ belongs to the ideal of compact
operators \cite{kl}, that is to $\tilde{J}_2$. For the inductive step, suppose that
$1-\sum_{j=1}^n z_jz_j^*$ is in $\tilde{J}_{2n}$. Then $1-\sum_{j=1}^{n+1} Z_jZ_j^*$
belongs to $\tilde{J}_{2n}\otimes\K$ by (\ref{oddradius}), and thus
to $\tilde{J}_{2(n+1)}$.
\end{proof}

It is natural to regard the quotient $C(B^{2n}_q)/J_{2n}$ as the
$C^*$-algebra of continuous functions on the `boundary' of the quantum ball
$B^{2n}_q$. Thus, we use the notation
\begin{equation}\label{b2nboundary}
C(\partial B^{2n}_q)=C(B^{2n}_q)/J_{2n}.
\end{equation}
Consequently, there is a natural surjection $\pi:C(B^{2n}_q)\ra
C(\partial B^{2n}_q)$ and we have an exact sequence
\begin{equation}\label{b2nexact}
0 \lra J_{2n} \lra C(B^{2n}_q) \overset{\pi}{\lra} C(\partial B^{2n}_q) \lra 0.
\end{equation}

Since $C(B^{2n}_q)\cong C^*(M_n)$ and ideal $J_{2n}$ is generated by
projection $P_{n+1}$, it follows from the general theory of
graph algebras that the quotient $C(\partial B^{2n}_q)$
is isomorphic to the $C^*$-algebra of a graph obtained from $M_n$ by
removing vertex $v_{n+1}$ and all edges $e_{i, n+1} \; ( 1\leq i\leq n)$.
In \cite[Section 4.1]{hs1}, this graph was denoted by $L_{2n-1}$.
For example, by removing vertex $v_4$ from graph $M_3$ we
obtain graph $L_5$, which looks as follows.

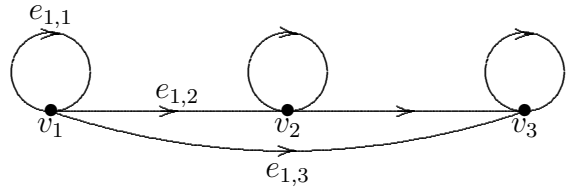
\[ \beginpicture
\setcoordinatesystem units <2.1cm,2.1cm>
\setplotarea x from -4 to 2, y from -0.7 to 0.6

\put {$\bullet$} at -3 0
\put {$\bullet$} at -1.5 0
\put {$\bullet$} at 0 0

\setlinear
\plot -3 0 -1.5 0 /
\plot -1.5 0 0 0 /

\arrow <0.235cm> [0.2,0.6] from -2.3 0 to -2.2 0
\arrow <0.235cm> [0.2,0.6] from -0.8 0 to -0.7 0
\setquadratic
\plot -3 0 -1.5 -0.25 0 0 /
\arrow <0.235cm> [0.2,0.6] from -1.5 -0.249 to -1.45 -0.24852
\circulararc 360 degrees from -3 0 center at -3 0.25
\circulararc 360 degrees from -1.5 0 center at -1.5 0.25
\circulararc 360 degrees from 0 0 center at 0 0.25
\arrow <0.235cm> [0.2,0.6] from -3 0.492 to -2.95 0.492
\arrow <0.235cm> [0.2,0.6] from -1.49 0.495 to -1.44 0.49
\arrow <0.235cm> [0.2,0.6] from 0 0.49 to 0.05 0.49
\put {$e_{1,1}$} [v] at -3 0.6
\put {$e_{1,2}$} [v] at -2.2 0.1
\put {$e_{1,3}$} [v] at -1.5 -0.4
\put {$v_1$} [v] at -3 -0.1
\put {$v_2$} [v] at -1.5 -0.1
\put {$v_3$} [v] at 0 -0.1
\endpicture \]

The Cuntz-Krieger generators of $C^*(L_{2n-1})$
are projections $\{Q_i:i=1,\ldots,n\}$ and partial isometries
$\{R_{i,j}:i=1,\ldots,n,\;j=i,\ldots,n\}$. Note that
$R_{n,n}$ is a partial unitary with domain and range projection $Q_n$.
It is worth mentioning that the above identification and
the results of \cite{hs1} immediately imply that
$C(\partial B^{2(n+1)}_q)$ is isomorphic with the quantum
double suspension of $C(\partial B^{2n}_q)$.

Let $\pi:C(B^{2n}_q)\ra C(\partial B^{2n}_q)=C(B^{2n}_q)/J_{2n}$ be
the natural surjection. We define a polynomial algebra of
$\partial B^{2n}_q$ as the $\pi$ image of the polynomials on $B^{2n}_q$,
that is ${\mathcal O}(\partial B^{2n}_q)=\pi({\mathcal O}(B^{2n}_q))$ (cf. Remark
\ref{algebraicNC}). We denote by $w_i=\pi(z_i)$, $i=1,\ldots,n$, the
generators of ${\mathcal O}(\partial B^{2n}_q)$. It follows from Theorem \ref{b2n}
and Lemma \ref{J2n} that $C^*$-algebra $C(\partial B^{2n}_q)$ has the
following presentation.

\begin{prop}\label{oddsphere}
$C(\partial B^{2n}_q)$ is isomorphic with the $C^*$-algebra
$C^*(w_1,\ldots,w_n)$, universal for the relations
\begin{eqnarray}
w_iw_j & = & q^{1/2}w_jw_i  \;\;\;\;\; \mbox{for} \;\; i<j,
\label{oddspherepresentation1} \\
w_iw_j^* & = & q^{-1/2}w_j^*w_i  \;\;\;\;\; \mbox{for} \;\; i<j,
\label{oddspherepresentation2} \\
w_i^*w_i-qw_iw_i^* & = & (1-q)\left(1-\sum_{j=i+1}^n w_jw_j^*\right)
\;\;\;\;\; \mbox{for} \;\; i=1,\ldots,n, \label{oddspherepresentation3} \\
\sum_{j=1}^n w_jw_j^* & = & 1. \label{oddspherepresentation4}
\end{eqnarray}
\end{prop}

Combining (\ref{oddspherepresentation3})
with (\ref{oddspherepresentation4}) one deduces
that element $w_1$ is normal and
\begin{equation}\label{oddspherepresentation5}
w_i^*w_i-w_iw_i^*=(1-q)\sum_{j=1}^{i-1}w_jw_j^*
\;\;\;\;\;\; \text{for} \;\; i=2,\ldots,n.
\end{equation}

It turns out that the boundary $\partial B^{2n}_q$ of
our noncommutative ball $B^{2n}_q$ is identical with the
unitary quantum sphere $S^{2n-1}_\mu$ of Vaksman-Soibelman \cite{vs}
both on the $C^*$-algebra and on the polynomial algebra level.
Indeed, in the convention of \cite{hs1}, ${\mathcal O}(S^{2n-1}_\mu)$ is
generated by $n$ elements $z_1,\ldots,z_n$ satisfying
relations (4.1)--(4.4) of \cite[Section 4]{hs1}.
We obtain an identification of this algebra with our
${\mathcal O}(\partial B^{2n}_q)$ by setting $\mu^2=q$, $z_j=w_{n-j+1}$ for
$j=1,\ldots,n-1$, and $z_n=w_1^*$. In particular,
${\mathcal O}(\partial B^{4}_{q^2})$ is identical with ${\mathcal O}(SU_q(2))$
of Woronowicz \cite{w}.

\subsection{The uniqueness criteria}

When working with algebras defined by universal properties it is not
difficult to construct their homomorphisms. However, it is usually a much
harder task to decide if a homomorphism is injective or not.
Very convenient criteria of injectivity of homomorphisms (known as
uniqueness theorems) have been developed for the class of graph algebras,
to which both $C(B^{2n}_q)$ and $C(\partial B^{2n}_q)$ belong. Thus,
the general uniqueness theorem for graph algebras \cite[Theorem 1.2]{sz}
implies injectivity criteria of homomorphisms of both $C(B^{2n}_q)$ and
$C(\partial B^{2n}_q)$. Propositions \ref{b2nuniqueness} and \ref{s2n1uniqueness}
contain reformulations of these criteria in terms of the generators provided
by Theorem \ref{b2n} and Proposition \ref{oddsphere} of the present paper.

\begin{prop}\label{b2nuniqueness}
A $*$-homomorphism $\alpha$ from $C(B^{2n}_q)$ into another $C^*$-algebra
is injective if and only if $\alpha(z_1)$ is not normal.
\end{prop}
\begin{proof}
Since the element $1-\sum_{j=1}^nz_jz_j^*$ generates an essential, simple
ideal of $C(B^{2n}_q)$ (namely $J_{2n}$), a $*$-homomorphism $\alpha$
from $C(B^{2n}_q)$ to another $C^*$-algebra is injective provided
$\alpha(1-\sum_{j=1}^nz_jz_j^*)\neq0$. However, it follows from relation
(\ref{b2npresentation3}) that condition $\alpha(1-\sum_{j=1}^nz_jz_j^*)\neq0$
is equivalent to the requirement that $\alpha(z_1)\alpha(z_1)^*\neq
\alpha(z_1)^*\alpha(z_1)$.
\end{proof}

\begin{prop}\label{s2n1uniqueness}
A $*$-homomorphism $\alpha$ from $C(\partial B^{2n}_q)$ into another
$C^*$-algebra is injective if and only if the spectrum of $\alpha(w_1)$
contains the entire unit circle.
\end{prop}
\begin{proof}
Let $\alpha$ be a $*$-homomorphism of $C(\partial B^{2n}_q)$. After the
identification of $C(\partial B^{2n}_q)$ with $C^*(L_{2n-1})$,
\cite[Theorem 1.2]{sz} implies that $\alpha$ is injective if and only
if the spectrum of $\alpha(R_{n,n})$ contains the entire unit circle.
By \cite[Theorem 4.4]{hs1}, this is equivalent to the requirement that
the spectrum of $\alpha(b_n)$ (with $b_n$ the normal generator of
$C(S^{2n-1}_\mu)=C^*(b_1,\ldots,b_n)$) contains the entire unit circle.
Since $b_n$ is normal and $b_n^*$ corresponds to $w_1$ under the identification
of $C(S^{2n-1}_\mu)$ with $C(\partial B^{2n}_q)$, the proposition follows.
\end{proof}


\section{The `odd-dimensional' noncommutative balls}

\subsection{The algebra}

Using the quantum double suspension, the $C^*$-algebras
$C(B^{2n-1}_q)$ of the `odd-dimensional' noncommutative balls
are defined inductively for all $n$ as follows:
$$ C(B^1_q)=C([-1,1]),\;\;\;\;\; C(B^{2n+1}_q)=\Sigma^2C(B^{2n-1}_q). $$

This definition, based on the general approach developed in \cite{hs1},
first appeared in \cite{f}. It is shown therein that the
primitive ideal space of $C(B^{2n-1}_q)$ consists of $n-1$ circles and a
closed interval, with certain natural non-Hausdorff topology. This fact
combined with the description of the primitive ideal spaces of all graph
$C^*$-algebras, given in \cite{hs2}, implies that $C(B^{2n-1}_q)$
are not isomorphic to any graph algebras. Thus, contrary to the
`even-dimensional' case, in the present situation we cannot rely on the
well-established machinery of graph algebras. However, some
structural properties of the $C^*$-algebras $C(B^{2n-1}_q)$ may be
deduced from their inductive definition via the quantum double suspension.
In particular, their $K$-groups can be calculated this way and shown to
coincide with those of the classical balls.

\begin{prop}\label{boddK}
We have $K_1(C(B^{2n-1}_q))=0$, and $K_0(C(B^{2n-1}_q))\cong\Z$
generated by the class of identity $[1]_0$.
\end{prop}
\begin{proof}
This is easily proved by induction on $n$, applying the six-term exact
sequence of $K$-theory to extension (\ref{doublesuspension}).
\end{proof}

Now we define generators of $C(B^{2n-1}_q)$ and thus introduce
their polynomial algebras. The $C^*$-algebra $C([-1,1])$ is generated by
one element $x_1$ satisfying $x_1=x_1^*$ and $x_1^2\leqq1$.
Applying inductively the method of Section 2, we obtain the generators
and relations for $C(B^{2n-1}_q)$, as follows.

\begin{theo}\label{bodd}
$C(B^{2n-1}_q)$ is isomorphic with the $C^*$-algebra
$C^*(x_1,\ldots,x_n)$, universal for the relations
\begin{eqnarray}
x_1 & = & x_1^*,
\label{boddpresentation1} \\
x_ix_j & = & q^{1/2}x_jx_i  \;\;\;\;\; \mbox{for} \;\; i<j,
\label{boddpresentation2} \\
x_ix_j^* & = & q^{-1/2}x_j^*x_i  \;\;\;\;\; \mbox{for} \;\; 2\leqq i<j,
\label{boddpresentation3} \\
x_i^*x_i-qx_ix_i^* & = & (1-q)\left(1-\sum_{j=i+1}^n x_jx_j^*\right)
\;\;\;\;\; \mbox{for} \;\; i=2,\ldots,n,
\label{boddpresentation4} \\
x_1^2+\sum_{j=2}^n x_jx_j^* & \leqq & 1.
\label{boddpresentation5}
\end{eqnarray}
\end{theo}
\begin{proof}
The following line of proof is very similar to the one from Theorem
\ref{b2n}. Therefore we only sketch the main points.

Contrary to relations (\ref{b2npresentation1})--(\ref{b2npresentation3})
of the `even-dimensional' noncommutative balls which admit representations
in unbounded operators, in the present case we always have $||x_j||\leq1$
due to (\ref{boddpresentation5}). To establish the required isomorphism,
we proceed by induction on $n$. Case $n=1$ is obvious. For the inductive step,
suppose $C(B^{2n-1}_q)\cong C^*(x_1,\ldots,x_n)$, and let $C^*(X_1,\ldots,
X_{n+1})$ be the universal $C^*$-algebra for relations
(\ref{boddpresentation1})--(\ref{boddpresentation5}), with $n$ replaced by $n+1$.
We must show that there exist $*$-homomorphisms
$\phi:C^*(X_1,\ldots,X_{n+1})\ra \Sigma^2C^*(x_1,\ldots,x_n)$ and
$\psi:\Sigma^2C^*(x_1,\ldots,x_n)\ra C^*(X_1,\ldots,X_{n+1})$ such
that $\psi\circ\phi=\id$ and $\phi\circ\psi=\id$. These two maps are
defined as in the proof of Theorem \ref{b2n}. Namely, $\phi$ maps $X_1,\ldots,
X_{n+1}$ to the $n+1$ generators of $\Sigma^2C^*(x_1,\ldots,x_n)$ given by formulae
(\ref{s2generators1})--(\ref{s2generators2}). Obviously, elements $\phi(X_1),
\ldots,\phi(X_{n+1})$ satisfy (\ref{boddpresentation1})--(\ref{boddpresentation4}).
They fulfill (\ref{boddpresentation5}) thanks to identity (\ref{oddradius}).
As in the proof of Theorem \ref{b2n}, one verifies that $X_{n+1}$ admits a
polar decomposition $T|X_{n+1}|$ in $C^*(X_1,\ldots,X_{n+1})$. Then a map
$\tilde{\psi}:C^*(x_1,\ldots,x_n)\ra C^*(X_1,\ldots,X_{n+1})$ is defined by
$\tilde{\psi}(x_i)=X_i(1-TT^*)$ for $i=1,\ldots,n$. Finally, $\psi$ is
constructed with help of Proposition \ref{universalqds} so that $\psi(1\otimes V)=T$ and
$\psi(a\otimes E_{00})=\tilde{\psi}(a)$ for all $a\in C^*(x_1,\ldots,x_n)$.
Verification of the identities $\phi\circ\psi=\id$ and $\psi\circ\phi=\id$
is carried out in the same way as in Theorem \ref{b2n}. In particular,
Lemma \ref{sumidentities} remains valid for operators $X_j$, $j=1,\ldots,n+1$,
satisfying (\ref{boddpresentation1})--(\ref{boddpresentation5}).
\end{proof}

In the remainder of this paper, we simply identify $C(B^{2n-1}_q)=
C^*(x_1,\ldots,x_n)$. We define the polynomial algebra
${\mathcal O}(B^{2n-1}_q)$ as the $*$-subalgebra of $C(B^{2n-1}_q)$
generated by $x_1,\ldots,x_n$.

\subsection{The boundary}

We define $J_{2n-1}$ as the closed, two-sided ideal of $C(B^{2n-1}_q)$
generated by $1-x_1^2-\sum_{j=2}^n x_jx_j^*$. In order to match the
requirements of Section 1, this ideal must be essential.

\begin{lemm}\label{Jodd}
$J_{2n-1}$ is an essential ideal of $C(B^{2n-1}_q)$. Furthermore,
$J_{2n+1}=J_{2n-1}\otimes\K$.
\end{lemm}
\begin{proof}
At first we observe that $J_{2n+1}=J_{2n-1}\otimes\K$ holds by
identity (\ref{oddradius}), since $x_1$ is self-adjoint.
Now we proceed by induction on $n$. If $n=1$
then $C(B^1_q)=C([-1,1])$ and $J_1=\{f\in C([-1,1]):f(-1)=f(1)=0\}$. Thus
$J_1$ is essential in $C(B^1_q)$. For the inductive step, suppose that
$J_{2n-1}$ is essential in $C(B^{2n-1}_q)$. Then $J_{2n+1}=J_{2n-1}\otimes\K$
is essential in $C(B^{2n-1}_q)\otimes\K$. Since, by definition of
the quantum double suspension, $C(B^{2n-1}_q)\otimes\K$ is an essential
ideal of $C(B^{2n+1}_q)=\Sigma^2C(B^{2n-1}_q)$, the inductive step follows.
\end{proof}
\noindent
As an immediate corollary of Lemma \ref{Jodd} we see that
$J_{2n-1}\cong C_0(\R)\otimes\K$ for $n\geq2$.

\vspace{2mm}
Analogously to the `even-dimensional' case we define the boundary of
$B^{2n-1}_q$ by taking quotient of its $C^*$-algebra with $J_{2n-1}$,
as follows.
$$ C(\partial B^{2n-1}_q)=C(B^{2n-1}_q)/J_{2n-1}. $$
Then, with the natural surjection $\pi:C(B^{2n-1}_q)\ra C(\partial B^{2n-1}_q)$,
we have an exact sequence
\begin{equation}\label{boddexact}
0 \lra J_{2n-1} \lra C(B^{2n-1}_q) \overset{\pi}{\lra}
C(\partial B^{2n-1}_q) \lra 0.
\end{equation}

We define the polynomial algebra of $\partial B^{2n-1}_q$ as the $\pi$
image of the polynomials on $B^{2n-1}_q$, that is ${\mathcal O}(\partial B^{2n-1}_q)=
\pi({\mathcal O}(B^{2n-1}_q))$. We denote by $t_i=\pi(x_i)$, $i=1,\ldots,n$, the
generators of ${\mathcal O}(\partial B^{2n-1}_q)$. In terms of these generators,
the $C^*$-algebra $C(\partial B^{2n-1}_q)$ has the following presentation,
which follows immediately from Theorem \ref{bodd} and our definitions.

\begin{prop}\label{evensphere}
$C(\partial B^{2n-1}_q)$ is isomorphic with the $C^*$-algebra
$C^*(t_1,\ldots,t_n)$, universal for the relations
\begin{eqnarray}
t_1 & = & t_1^*, \label{evenspherepresentation1} \\
t_it_j & = & q^{1/2}t_jt_i  \;\;\;\;\; \mbox{for} \;\; i<j,
\label{evenspherepresentation2} \\
t_it_j^* & = & q^{-1/2}t_j^*t_i  \;\;\;\;\; \mbox{for} \;\; 2\leqq i<j,
\label{evenspherenpresentation3} \\
t_i^*t_i-qt_it_i^* & = & (1-q)\left(1-\sum_{j=i+1}^n t_jt_j^*\right)
\;\;\;\;\; \mbox{for} \;\; i=2,\ldots,n, \label{evenspherepresentation4} \\
t_1^2+\sum_{j=2}^n t_jt_j^* & = & 1. \label{evenspherepresentation5}
\end{eqnarray}
\end{prop}

It is useful to observe that combining (\ref{evenspherepresentation4})
with (\ref{evenspherepresentation5}) one obtains
\begin{equation}\label{oddspherepresentation6}
t_i^*t_i-t_it_i^*=(1-q)\left(t_1^2+\sum_{j=2}^{i-1}t_jt_j^*\right).
\end{equation}

In the case of $n=2$ we have the following presentation of
$C(\partial B^3_q)=C^*(t_1,t_2)$:
$$ t_1=t_1^*, \;\;\; t_1t_2=q^{1/2}t_2t_1, \;\;\;
   t_2^*t_2-qt_2t_2^*=1-q, \;\;\; t_1^2+t_2t_2^*=1. $$
These are exactly the relations defining the equatorial Podle\'s
sphere \cite[Formulae (7b)]{p}, with identification
$t_1=A$, $t_2=B^*$, and $q=\mu^4$. More generally, our algebras
${\mathcal O}(\partial B^{2n-1}_q)$ are easily seen to coincide with the
algebras $A(S^{2n-2}_\mu)$ of `even-dimensional' Euclidean quantum
spheres studied in \cite{frt,hl}. Indeed, using the presentation
for $A(S^{2n-2}_\mu)$ in terms of generators $c_0,\ldots,c_{n-1}$
(with quantization parameter $\mu$) given in \cite[Section 2]{hl},
we obtain an identification of this algebra with our
${\mathcal O}(\partial B^{2n-1}_q)$ by setting $\mu^2=q$ and $c_{j-1}=\mu^{n-j}t_j$
for $j=1,\ldots,n$.

\begin{rema}\label{evenspheresgraph}
\rm It is worth noting that unlike the $C^*$-algebras $C(B^{2n-1}_q)$ of
the noncommutative `odd-dimensional' balls, the $C^*$-algebras $C(\partial B^{2n-1}_q)$
of their boundaries are isomorphic with certain graph algebras. Indeed, it follows immediately
from Propositions \ref{oddsphere} and \ref{evensphere} that $C(\partial B^{2n-1}_q)$
is isomorphic to the quotient of $C(\partial B^{2n}_q)$ by the ideal generated by
$w_1-w_1^*$. By virtue of our identification of $C(\partial B^{2n}_q)$ with
$C(S^{2n-1}_mu)$ and \cite[Proposition 5.1]{hs1}, we have $C(\partial B^{2n-1}_q)
\cong C^*(L_{2n-2})$ where $L_{2n-2}$ is the directed graph described therein.
This fact and \cite[Example 6.4]{hs1} imply that
$C^*$-algebra $C(\partial B^{2n+1}_q)$ is isomorphic
with the quantum double suspension $\Sigma^2 C(\partial B^{2n-1}_q)$.
\end{rema}

Since $C(\partial B^{2n-1}_q)$ is isomorphic to a graph algebra,
\cite[Theorem 1.2]{sz} yields a criterion of injectivity of its
homomorphisms. The following is a reformulation of this criterion
in terms of the generators from Proposition \ref{evensphere}.

\begin{prop}\label{s2nuniqueness}
A $*$-homomorphism $\alpha$ from $C(\partial B^{2n-1}_q)$ into another
$C^*$-algebra is injective if and only if the spectrum of $\alpha(t_1)$
contains both positive and negative numbers.
\end{prop}


\section{The `even-dimensional' glued quantum spheres}

For any complex numbers $\lambda_1,\ldots,\lambda_n$ of modulus one there
exists a $*$-automorphism $\beta$ of ${\mathcal O}(\partial B^{2n}_q)$ such that
\begin{equation}\label{betaodd-1}
\beta(w_j)=\lambda_jw_j \;\;  \mbox{for} \; j=1,\ldots,n, \;\; \text{or}
\end{equation}
\begin{equation}\label{betaodd-2}
\beta(w_1)=\lambda_1w_1^*, \; \; \;
\beta(w_j)=\lambda_jw_j \;\; \mbox{for}\; j=2,\ldots,n .
\end{equation}
Any such an automorphism extends to the $C^*$-algebra $C(\partial B^{2n}_q)$
and the extension is still denoted $\beta$.
We define $C^*$-algebras $C(S^{2n}_{q,\beta})$ of the `even-dimensional'
glued quantum spheres as the corresponding doubles
$$ C(S^{2n}_{q,\beta})=C(B^{2n}_q)\oplus_\beta C(B^{2n}_q), $$
according to the general recipe given in (\ref{double}).
These are type $I$ $C^*$-algebras and the exact sequence (\ref{doubleexact})
in the present case takes the form
\begin{equation}\label{s2nexact}
0 \lra J_{2n}\oplus J_{2n} \lra C(S^{2n}_{q,\beta}) \lra C(\partial B^{2n}_q) \lra 0,
\end{equation}
with $J_{2n}\cong\K$. Consequently, regardless of the choice of $\beta$,
the primitive ideal space of $C(S^{2n}_{q,\beta})$ consists of two points and
$n$ circles, with certain non-Hausdorff topology.

If $\beta_1$ and $\beta_2$ are automorphisms of $C(\partial B^{2n}_q)$
both of type (\ref{betaodd-1}) or both of type (\ref{betaodd-2}),
respectively, then $\beta_1^{-1}\circ\beta_2$ is of the form
(\ref{betaodd-1}). Since any automorphism of type (\ref{betaodd-1})
admits a lift to an automorphism  of $C(B^{2n}_q)$,
Proposition \ref{doubletype} implies that the choice of
the scalars $\lambda_j$ does not affect the isomorphism class
of $C(S^{2n}_{q,\beta})$. Consequently, it suffices to consider
two cases only: $\beta=\id$ and $\beta(w_1)=w_1^*$, $\beta(w_j)=
w_j$ for $j=2,\ldots,n$. We will show, below, that these two
choices yield non-isomorphic $C^*$-algebras. While $\beta=\id$
gives rise to previously known quantum spheres, the latter case
produces a new class of quantum spheres which we call {\em mirror
quantum spheres}. This construction and analysis generalizes results
from \cite{hms2}, applicable to the case of $C(S^2_{q,\beta})$.

We define the polynomial algebra ${\mathcal O}(S^{2n}_{q,\beta})$ as follows. If
$\beta$ is of the form (\ref{betaodd-1}), then ${\mathcal O}(S^{2n}_{q,\beta})$
is the $*$-subalgebra of $C(S^{2n}_{q,\beta})$ generated by
\begin{eqnarray}
e_0 & = & \left(\left(1-\sum_{j=1}^n z_jz_j^*\right)^{1/2},\; \;
-\left(1-\sum_{j=1}^n z_jz_j^*\right)^{1/2}\right),
\label{s2ngenerators2} \\
e_i & = & (\lambda_iz_i,z_i), \;\; \mbox{for} \; i=1,\ldots,n.
\label{s2ngenerators1}
\end{eqnarray}
Note that element $e_0$ in (\ref{s2ngenerators2}) is well-defined by virtue of inequality
(\ref{evenradiusinequality}). If $\beta$ is of the form (\ref{betaodd-2}),
then ${\mathcal O}(S^{2n}_{q,\beta})$ is the $*$-subalgebra of $C(S^{2n}_{q,\beta})$
generated by
\begin{eqnarray}
e_0 & = & \left(\left(1-z_1^*z_1-\sum_{j=2}^n z_jz_j^*\right)^{1/2},\; \;
-\left(1-\sum_{j=1}^n z_jz_j^*\right)^{1/2}\right), \label{s2ngenerators5} \\
e_1 & = & (\lambda_1z_1^*,z_1), \label{s2ngenerators3} \\
e_i & = & (\lambda_iz_i,z_i), \;\; \mbox{for} \; i=2,\ldots,n. \label{s2ngenerators4}
\end{eqnarray}
Again, note that element $e_0$ in (\ref{s2ngenerators5}) is well-defined, since
$1-z_1^*z_1-\sum_{j=2}^n z_jz_j^*\geq 0$,
or equivalently $z_1^*z_1\geq z_1z_1^*$ (use (\ref{evenradiusinequality})
and (\ref{b2npresentation3})). This last inequality is proved
for $C(B^{2n}_q)$ by induction on $n$. Indeed,
$1-z_1^*z_1=q(1-z_1z_1^*)\leqq 1-z_1z_1^*$ in $C(B^2_q)$, and for
the inductive step use (\ref{b2nphi1}).

In either case, it is not difficult to verify that ${\mathcal O}(S^{2n}_{q,\beta})$
is a dense $*$-subalgebra of $C(S^{2n}_{q,\beta})$.

\vspace{2mm}
We now show that automorphisms of type (\ref{betaodd-1}) lead to quantum
spheres $S^{2n}_{q,\beta}$ identical with the previously discussed
boundaries of noncommutative `odd-dimensional' balls. To this end,
we prove that there exists a $C^*$-algebra isomorphism from
$C(\partial B^{2n+1}_q)$ to $C(S^{2n}_{q,\beta})$ which preserves
their polynomial algebras.

\begin{prop}\label{s2nid}
Let $\beta$ be a $*$-automorphism of $C(\partial B^{2n}_q)$ of the
form (\ref{betaodd-1}). Then there exists an isomorphism
$$ \phi:C(\partial B^{2n+1}_q)\ra C(S^{2n}_{q,\beta}) $$
such that $\phi(t_i)=e_{i-1}$ for $i=1,\ldots,n+1$.
\end{prop}
\begin{proof}
Such a $*$-homomorphism $\phi:C(\partial B^{2n+1}_q)\ra
C(S^{2n}_{q,\beta})$ exists by the universal property of
$C(\partial B^{2n+1}_q)$ from Proposition \ref{evensphere}. Indeed,
elements $e_0,e_1,\ldots,e_n$ of $C(S^{2n}_{q,\beta})$ satisfy
relations (\ref{evenspherepresentation1})--(\ref{evenspherepresentation5})
for $t_1,\ldots,t_{n+1}$. The only non-trivial condition $e_0e_i=q^{1/2}e_ie_0$,
$i=1,\ldots n$, holds due to the following identity
satisfied in $C(B^{2n}_q)=C^*(z_1, z_2, \cdots , z_n)$:
\begin{equation}\label{e0com}
\left(1-\sum_{j=1}^n z_jz_j^*\right)^{1/2}z_i=q^{1/2}z_i
\left(1-\sum_{j=1}^n z_jz_j^*\right)^{1/2} \;\;\; \mbox{for all}
\;\; i=1,\ldots,n,
\end{equation}
which may be verified by a straightforward induction based
on (\ref{oddradius}).

Surjectivity of $\phi$ is obvious, while its injectivity follows
from Proposition \ref{s2nuniqueness}.
\end{proof}

\begin{rema}\label{evencomp}
\rm In view of Proposition \ref{s2nid} and the discussion following
Proposition \ref{evensphere}, we may conclude that the same
`even-dimensional' quantum spheres can be obtained through
one of the following four distinct constructions:
\begin{description}
\item[(i)] as homogeneous spaces of the quantum orthogonal groups,
\item[(ii)] as boundaries of `odd-dimensional' noncommutative balls,
\item[(iii)] by gluing `even-dimensional' noncommutative balls along
their boundaries, and
\item[(iv)] by repeated application of the quantum double suspension
applied to the classical $2$-point space.
\end{description}
\end{rema}

If the boundary automorphism is of the form (\ref{betaodd-2}),
then it is not clear whether the generators of $C(S^{2n}_{q,\beta})$
have universal property with respect to a finite set of algebraic relations,
as is the case with $C(S^{2n}_{q,\id})\cong C(\partial B^{2n+1}_q)$
(see Propositions \ref{s2nid} and \ref{evensphere}). In the special
case of $C(S^2_{q,\beta})$, such a presentation was given in \cite{hms2}
after enlarging the polynomial algebra by the positive $e_0^+$ and the
negative $e_0^-$ parts of $e_0$. However, elements $e_0^+,e_0^-$
correspond to continuous but not differentiable functions. Thus we
prefer not to include them in our polynomial algebra ${\mathcal O}(S^{2n}_{q,\beta})$.

\vspace{2mm}
Our next goal is comparison of the $C^*$-algebras corresponding
to the two distinct forms of the boundary automorphism $\beta$.
This will be achieved by a careful analysis of their $K$-theory.
A natural basis for calculating the $K$-theory of $C(S^{2n}_{q,\beta})
=C(B^{2n}_q)\oplus_\beta C(B^{2n}_q)$ is the exact sequence
(\ref{s2nexact}).

\begin{theo}\label{s2nbetaK}
Let $\beta$ be a $*$-automorphism of $C(\partial B^{2n}_q)$ of the
form {\em (\ref{betaodd-1})} or {\em (\ref{betaodd-2})}. We have
$$ K_0(C(S^{2n}_{q,\beta}))\cong\Z^2 \;\;\; \text{and} \;\;\;
   K_1(C(S^{2n}_{q,\beta}))=0. $$
Generators of the $K_0$ group depend on $\beta$ as follows.
\begin{description}
\item[(i)] If $\beta$ is of the form {\em (\ref{betaodd-1})} then the
$K_0$ group is generated by $[1]_0$ and $[p_1]_0=-[p_2]_0$.
\item[(ii)] If $\beta$ is of the form {\em (\ref{betaodd-2})} then the $K_0$
group is generated by $[1]_0$ and $[p_1]_0=[p_2]_0$.
\end{description}
Herein, $p_1$ and $p_2$ are minimal projections in $J_{2n}\oplus 0$ and
$0\oplus J_{2n}$, respectively.
\end{theo}
\begin{proof}
If $\beta$ is of the form (\ref{betaodd-1}), then these claims follow
from the identifications summarized in Remark \ref{evencomp} and
the results of \cite{hs1}. Thus, we may consider only the case when
$\beta$ is of the form (\ref{betaodd-2}). Furthermore, we may take
all the scalars $\lambda_i$ to be equal to $1$, and thus
$\beta(w_1)=w_1^*$, $\beta(w_j)=w_j$ for $j=2,\ldots,n$.

It is more convenient to view $C(B^{2n}_q)$
as the graph algebra $C^*(M_n)$, and thus
$C(S^{2n}_{q,\beta})=C^*(M_n)\oplus_\beta C^*(M_n)$. Likewise
$C(\partial B^{2n}_q)$ is the graph algebra $C^*(L_{2n-1})$, as
explained in Section 3. Then in terms of the Cuntz-Krieger generators of
$C^*(L_{2n-1})$, the automorphism $\beta$ acts as $\beta(R_{i,j})=
R_{i,j}$ for $i=1,\ldots,n-1$, $j=i,\ldots,n$, and $\beta(R_{n,n})=R^*_{n,n}$.
This is easily seen by combining \cite[Theorem 4.4]{hs1} and
our identification of $C(\partial B^{2n}_q)$ with $C(S^{2n-1}_\mu)$
given below Proposition \ref{oddsphere}. Now sequence (\ref{s2nexact})
takes the form
\begin{equation}\label{s2ngraphexact}
0 \lra \K\oplus\K \lra C^*(M_n)\oplus_\beta C^*(M_n) \lra C^*(L_{2n-1}) \lra 0.
\end{equation}
Applying the six-term exact sequence of $K$-theory we get
\begin{equation}\label{6term}
\begin{CD}
K_0(\K\oplus\K) \cong \Z^2 @>  >> K_0(C(S^{2n}_{q,\beta})) @> >> \Z \\
@A {\partial_{\ind}} AA @. @VV  V \\
K_1(C^*(L_{2n-1})) \cong \Z @<<  < K_1(C(S^{2n}_{q,\beta})) @<< < 0
\end{CD}
\end{equation}
We must determine the index map. To this end, note that the
$K_1$ group of $C^*(L_{2n-1})$ is generated by the class of unitary
$\tilde{U}=R_{n,n}+1-Q_n$, due to R{\o}rdam's description of generators
of the $K_1$ group of a Cuntz-Krieger algebra \cite{ro}. This $\tilde{U}$
lifts to a partial isometry $U=(S_{n,n}+1-P_n,S^*_{n,n}+1-P_n)$ in
$C^*(M_n)\oplus_\beta C^*(M_n)$. Hence
$$ \partial_{\ind}([\tilde{U}]_1)=[1-U^*U]_0-[1-UU^*]_0=
   [(0,P_{n+1})]_0-[(P_{n+1},0)]_0\in K_0(\K\oplus\K).  $$
Since $P_{n+1}$ is a minimal projection in the ideal $\K$, we have
$\partial_{\ind}(1)=(-1,1)$. This immediately implies that
$K_0(C(S^{2n}_{q,\beta}))\cong\Z^2$ and $K_1(C(S^{2n}_{q,\beta}))=0$.

Since $K_0(C(\partial B^{2n}_q))$ is generated by $[1]_0$,
$K_0(C(S^{2n}_{q,\beta}))$ is generated by the classes of minimal
projections in $\K\oplus\K$ and the class of identity. However, in
$K_0(C^*(M_n)\oplus_\beta C^*(M_n))$ we have
$$ 0=[1-U^*U]_0-[1-UU^*]_0=[(0,P_{n+1})]_0-[(P_{n+1},0)]_0, $$
and hence $[(0,P_{n+1})]_0=[(P_{n+1},0)]_0$. Consequently,
$K_0(C(S^{2n}_{q,\beta}))$ is generated by $[1]_0$ and $[p_1]_0=[p_2]_0$,
where $p_1$ and $p_2$ are minimal projections in $J_{2n}\oplus 0$ and
$0\oplus J_{2n}$, respectively.
\end{proof}
\noindent
Alternatively, the $K$-theory of $K_0(C(S^{2n}_{q,\beta}))$ may be
determined with help of the Mayer-Vietoris argument (cf. \cite{bhms}).

\begin{coro}\label{s2nbetaid}
Let $\beta$ be a $*$-automorphism of $C(\partial B^{2n}_q)$ of the
form {\em (\ref{betaodd-2})}. Then $C(S^{2n}_{q,\beta})$ is not stably
isomorphic to $C(S^{2n}_{q,\id})$.
\end{coro}

Comparing the $C^*$-algebras $C(S^{2n}_{q,\beta})$ with the boundary
automorphisms $\beta$ of the form (\ref{betaodd-1}) or (\ref{betaodd-2}),
we see that they are type $I$, have homeomorphic primitive ideal
spaces and isomorphic $K$-groups. And yet these $C^*$-algebras are
non-isomorphic. In addition, Corollary
\ref{s2nbetaid} combined with Proposition \ref{doubletype}
implies that the automorphism $\beta$ of $C(\partial B^{2n}_q)$ such that
$\beta(w_1)=w_1^*$, $\beta(w_i)=w_i$ for $i=2,\ldots, n$, does not
admit a lift to an automorphism of $C(B^{2n}_q)$.

\begin{coro}\label{s2nbetaH}
Let $\beta$ be a $*$-automorphism of $C(\partial B^{2n}_q)$ of the
form {\em (\ref{betaodd-1})} or {\em (\ref{betaodd-2})}. Then the
$K$-homology groups are $K^0(C(S^{2n}_{q,\beta}))\cong\Z^2$ and
$K^1(C(S^{2n}_{q,\beta}))=0$.
\end{coro}
\begin{proof}
Since $K_0(C(S^{2n}_{q,\beta}))\cong\Z^2$ and $K_1(C(S^{2n}_{q,\beta}))=0$,
the Universal Coefficient Theorem \cite{rs} immediately implies that also
$K^0(C(S^{2n}_{q,\beta}))\cong\Z^2$ and $K^1(C(S^{2n}_{q,\beta}))=0$.
\end{proof}


\section{The `odd-dimensional' glued quantum spheres}

We now briefly go over the case of `odd-dimensional' glued quantum spheres.
Unlike the previously discussed `even-dimensional' case, this time
we do not obtain any new examples of quantum spheres.

For any complex numbers $\lambda_2,\ldots,\lambda_n$ of modulus one there
exists a $*$-automorphism $\beta$ of ${\mathcal O}(\partial B^{2n-1}_q)$ such that
\begin{equation}\label{betaeven-1}
\beta(t_1)=t_1, \; \; \beta(t_j)=\lambda_jt_j \;\;
\mbox{for} \; j=2,\ldots,n, \;\; \text{or}
\end{equation}
\begin{equation}\label{betaeven-2}
\beta(t_1)=-t_1, \; \; \beta(t_j)=\lambda_jt_j \;\;  \mbox{for} \; j=2,\ldots,n. \;\;
\end{equation}
Both automorphisms extend to the $C^*$-algebra $C(\partial B^{2n-1}_q)$
and the extensions are still denoted $\beta$. Let
$$ C(S^{2n-1}_{q,\beta})=C(B^{2n-1}_q)\oplus_\beta C(B^{2n-1}_q) $$
be the corresponding double, as defined in (\ref{double}).
Regardless of the choice of such an automorphism $\beta$, all of
the $C^*$-algebras $C(S^{2n-1}_{q,\beta})$ are isomorphic by
Proposition \ref{doubletype}.

We define the polynomial algebra ${\mathcal O}(S^{2n-1}_{q,\beta})$ as follows. If
$\beta$ is of the form (\ref{betaeven-1}), then ${\mathcal O}(S^{2n-1}_{q,\beta})$
is the $*$-subalgebra of $C(S^{2n-1}_{q,\beta})$ generated by
\begin{eqnarray}
f_0 & = & \left(\left(1-x_1^2-\sum_{j=2}^n x_jx_j^*\right)^{1/2}, \; \;
-\left(1-x_1^2-\sum_{j=2}^n x_jx_j^*\right)^{1/2}\right), \label{s2n1generators3} \\
f_1 & = & (x_1,x_1), \label{s2n1generators1} \\
f_i & = & (\lambda_ix_i,x_i), \;\; \mbox{for} \; i=2,\ldots,n. \label{s2n1generators2}
\end{eqnarray}
If $\beta$ is of the form (\ref{betaeven-2}), then
${\mathcal O}(S^{2n-1}_{q,\beta})$ is the analogous $*$-subalgebra of
$C(S^{2n-1}_{q,\beta})$, except
\begin{equation}\label{s2n1generators4}
f_1=(-x_1,x_1).
\end{equation}
In either case, ${\mathcal O}(S^{2n-1}_{q,\beta})$ is a dense $*$-subalgebra of
$C(S^{2n-1}_{q,\beta})$. Furthermore, it turns out that these algebras are
isomorphic with the polynomial algebras ${\mathcal O}(\partial B^{2n}_q)$ of the boundaries
of `even-dimensional' noncommutative balls. To prove this fact,
we show that there exists a $C^*$-algebra isomorphism from $C(\partial B^{2n}_q)$
to $C(S^{2n-1}_{q,\beta})$ which preserves their polynomial algebras.

\begin{prop}\label{s2n1}
Let $\beta$ be a $*$-automorphism of $C(\partial B^{2n-1}_q)$
of the form {\em (\ref{betaeven-1})} or {\em (\ref{betaeven-2})}.
Then there exists an isomorphism
$$ \phi:C(\partial B^{2n}_q)\ra C(S^{2n-1}_{q,\beta}) $$
such that $\phi(w_1)=f_0+if_1$ and $\phi(w_k)=f_{k}$ for $k=2,\ldots,n$.
\end{prop}
\begin{proof}
Such a $*$-homomorphism $\phi:C(\partial B^{2n}_q)\ra
C(S^{2n-1}_{q,\beta})$ exists by the universal property of
$C(\partial B^{2n}_q)$ from Proposition \ref{oddsphere}. Indeed,
elements $f_0+if_1,f_2,\ldots,f_n$ of $C(S^{2n-1}_{q,\beta})$ satisfy
relations (\ref{oddspherepresentation1})--(\ref{oddspherepresentation4})
for $w_1,\ldots,w_n$. This follows easily from our definitions and the
following two identities which hold in $C(B^{2n-1}_q)$:
\begin{eqnarray}
\left(1-x_1^2-\sum_{j=2}^n x_jx_j^*\right)^{1/2}x_1 & = &
x_1\left(1-x_1^2-\sum_{j=2}^n x_jx_j^*\right)^{1/2}, \label{f0com1} \\
\left(1-x_1^2-\sum_{j=2}^n x_jx_j^*\right)^{1/2}x_i
& = & q^{1/2}x_i\left(1-x_1^2-\sum_{j=2}^n x_jx_j^*\right)^{1/2},
\;\;\; i=2,\ldots,n. \label{f0com2}
\end{eqnarray}
Surjectivity of $\phi$ is clear, while its injectivity follows from Proposition
\ref{s2n1uniqueness}.
\end{proof}

\begin{rema}\label{oddcomp}
\rm In view of Proposition \ref{s2n1} and the discussion following
Proposition \ref{oddsphere}, we may conclude that the same
`odd-dimensional' quantum spheres can be obtained through one of the following
four distinct constructions:
\begin{description}
\item[(i)] as homogeneous spaces of the quantum unitary groups,
\item[(ii)] as boundaries of `even-dimensional' noncommutative balls,
\item[(iii)] by gluing `odd-dimensional' noncommutative balls along
their boundaries, and
\item[(iv)] by repeated application of the quantum double suspension
applied to the classical circle $S^1$.
\end{description}
\end{rema}


\section{Irreducible representations}

In this section, we give explicit formulae for irreducible
representations of the noncommutative balls $C(B^n_q)$ and the
mirror quantum spheres $C(S^{2n}_{q,\beta})$. Since the
$C^*$-algebras $C(\partial B^n_q)$ and $C(S^{2n-1}_{q,\beta})$
are isomorphic to the well-studied quantum spheres, their
irreducible representations are already available in the literature.

\subsection{Irreducible representations of the noncommutative balls}

Routine proofs of the following Propositions \ref{bevenirr} and
\ref{boddirr} are omitted. They are
established by induction based on the exact sequence
$$ 0 \lra C(B^n_q)\otimes\K \lra C(B^{n+2}_q) \lra C(S^1) \lra 0 $$
and formulae (\ref{s2generators1})--(\ref{s2generators2}). We use the following
notation. $\H_m$ denotes a Hilbert space with an orthonormal basis
$\{\xi_{k_1,\ldots,k_m}:k_i=0,1,\ldots\}$, and $S_r$ is a weighted shift on
$\H_m$ defined by
$$ S_r(\xi_{k_1,\ldots,k_m})= \begin{cases}
   \sqrt{(1-q^{1+k_r})q^{k_{r+1}+k_{r+2}\cdots+k_m}}
   \xi_{k_1,\ldots,k_{r-1},1+k_r,k_{r+1},\ldots,k_m}, & 1\leq r\leq m-1 \\
   \sqrt{1-q^{1+k_m}}\xi_{k_1,\ldots,k_{m-1},1+k_m}, & r=m. \end{cases} $$

\begin{prop}\label{bevenirr}
The following is a complete {\em (}up to unitary equivalence{\em )} list
of irreducible representations of the $C^*$-algebra $C(B^{2n}_q)=
C^*(z_1,\ldots,z_n)$:
$$ \irr(B^{2n}_q)=\{\1,\2,\3,\ldots,\n:\theta\in{\C},
   |\theta|=1\}\cup\{\sigma\}. $$
Herein representations $\1$ are $1$-dimensional, $\rho^\theta_j$,
$j=2,\ldots,n$, act on $\H_{j-1}$, and $\sigma$ acts on $\H_n$.
On the generators $z_1,\ldots,z_n$, these representations are
given by the following formulae:
\begin{eqnarray*}
\1(z_i) & = &
  \begin{cases}
    0, & 1\leq i\leq n-1 \\
    \theta, & i=n
  \end{cases} \\
\rho_j^\theta(z_i)\xi_{k_1,\ldots,k_{j-1}} & = &
  \begin{cases}
    0, & 1\leq i\leq n-j \\
    \theta q^{(k_1+k_2+\cdots+k_{j-1})/2}\xi_{k_1,\ldots,k_{j-1}}, & i=n-j+1 \\
    S_{i+j-n-1}(\xi_{k_1,\ldots,k_{j-1}}), & n-j+2\leq i\leq n
  \end{cases} \\
\sigma(z_i)\xi_{k_1,\ldots,k_n} & = & S_i(\xi_{k_1,\ldots,k_n}), \;\;\; 1\leq i\leq n
\end{eqnarray*}
\end{prop}
Considering exact sequence (\ref{b2nexact}), we see that
$\sigma$ of Proposition \ref{bevenirr} is an extension
to $C(B^{2n}_q)$ of the irreducible representation of the
ideal $J_{2n}\cong\K$, while $\rho_j^\theta$, $j=1,\ldots,n$
are the lifts of the irreducible representations of the
quotient $C(\partial B^{2n}_q)$.

\begin{prop}\label{boddirr}
The following is a complete {\em (}up to unitary equivalence{\em )} list
of irreducible representations of the $C^*$-algebra $C(B^{2n+1}_q)=
C^*(x_1,\ldots,x_{n+1})$:
$$ \irr(B^{2n+1}_q)=\{\5,\6,\7,\ldots,\eta_n^\theta:\theta\in{\C},
   |\theta|=1\}\cup\{\sigma_s:s\in[-1,1]\}. $$
Herein representations $\5$ are $1$-dimensional, $\eta^\theta_j$,
$j=2,\ldots,n$, act on $\H_{j-1}$, and $\sigma_s$ act on $\H_n$.
On the generators $x_1,\ldots,x_{n+1}$, these representations are
given by the following formulae:
\begin{eqnarray*}
\5(x_i) & = &
  \begin{cases}
    0, & 1\leq i\leq n \\
    \theta, & i=n+1
  \end{cases} \\
\eta_j^\theta(x_i)\xi_{k_1,\ldots,k_{j-1}} & = &
  \begin{cases}
    0, & 1\leq i\leq n-j+1 \\
    \theta q^{(k_1+k_2+\cdots+k_{j-1})/2}\xi_{k_1,\ldots,k_{j-1}}, & i=n-j+2 \\
    S_{i+j-n-2}(\xi_{k_1,\ldots,k_{j-1}}), & n-j+3\leq i\leq n+1
  \end{cases} \\
\sigma_s(x_i)\xi_{k_1,\ldots,k_n} & = &
  \begin{cases}
    sq^{(k_1+k_2+\cdots+k_n)/2}\xi_{k_1,\ldots,k_n}, & i=1 \\
    S_{i-1}(\xi_{k_1,\ldots,k_n}), & 2\leq i\leq n+1
  \end{cases}
\end{eqnarray*}
\end{prop}
Considering exact sequence (\ref{boddexact}), we see that
$\sigma_s$, $s\in(-1,1)$ of Proposition \ref{boddirr} are extensions
to $C(B^{2n+1}_q)$ of irreducible representations of the
ideal $J_{2n+1}\cong C_0(\R)\otimes\K$, while $\sigma_{\pm1}$ and
$\rho_j^\theta$, $j=1,\ldots,n$ are the lifts of irreducible
representations of the quotient $C(\partial B^{2n+1}_q)$.

\begin{rema}\label{newproof}
\rm A more direct proof of Proposition \ref{s2nuniqueness} may be
constructed with help of Proposition \ref{boddirr}, as follows.
Of all the irreducible representations of $C(B^{2n-1}_q)$ only
$\sigma_{\pm1}$ and $\rho_j^\theta$, $j=1,\ldots,n$
descend to $C(\partial B^{2n-1}_q)$. Then it is not difficult to observe that
$\sigma_1\oplus\sigma_{-1}$ is faithful on $C(\partial B^{2n-1}_q)$.
This immediately implies Proposition \ref{s2nuniqueness}.
\end{rema}

\subsection{Irreducible representations of the mirror quantum spheres}

We now turn to description of irreducible representations of the mirror
quantum spheres $C(S_{q,\beta}^{2n})$, corresponding to automorphisms
$\beta$ given by formula (\ref{betaodd-2}). It suffices to consider
the case with $\lambda_j=1$ for all $j=1,\ldots,n$.

\begin{prop}\label{mirrorirr}
Let $\beta$ be an automorphism of $C(\partial B^{2n}_q)$ given by
(\ref{betaodd-2}) with $\lambda_j=1$ for all $j=1,\ldots,n$.
The following is a complete {\em (}up to unitary equivalence{\em )} list of
irreducible representations of the $C^*$-algebra $C(S^{2n}_{q,\beta})=
C^*(e_0,e_1,\ldots,e_n)$:
$$ \irr(S^{2n}_{q,\beta})=\{\varrho_1^\theta,\varrho_2^\theta,
   \varrho_3^\theta,\ldots,\varrho_n^\theta:\theta\in{\C},
   |\theta|=1\}\cup\{\sigma_+,\sigma_-\}. $$
Herein representations $\varrho_1^\theta$ are $1$-dimensional, $\varrho^\theta_j$,
$j=2,\ldots,n$, act on $\H_{j-1}$, and $\sigma_+$, $\sigma_-$ act on $\H_n$.
On the generators $e_0,e_1,\ldots,e_n$, these representations are
given by the following formulae:
\begin{eqnarray*}
\varrho_1^\theta(e_i) & = &
  \begin{cases}
    0, & 0\leq i\leq n-1 \\
    \theta, & i=n
  \end{cases} \\
\text{for} \; 2\leq j\leq n-1, \;\; \varrho_j^\theta(e_i)\xi_{k_1,\ldots,k_{j-1}} & = &
  \begin{cases}
    0, & 0\leq i\leq n-j \\
    \theta q^{(k_1+k_2+\cdots+k_{j-1})/2}\xi_{k_1,\ldots,k_{j-1}}, & i=n-j+1 \\
    S_{i+j-n-1}(\xi_{k_1,\ldots,k_{j-1}}), & n-j+2\leq i\leq n
  \end{cases} \\
\varrho_n^\theta(e_i)\xi_{k_1,\ldots,k_{n-1}} & = &
  \begin{cases}
    0, & i=0 \\
    \overline{\theta}q^{(k_1+k_2+\cdots+k_{n-1})/2}\xi_{k_1,\ldots,k_{n-1}}, & i=1 \\
    S_{i-1}(\xi_{k_1,\ldots,k_{n-1}}), & 2\leq i\leq n
  \end{cases} \\
\sigma_+(e_i)\xi_{k_1,\ldots,k_n} & = & \begin{cases}
  q^{(1+k_1+k_2+\cdots+k_n)/2}\xi_{k_1,\ldots,k_n}, & i=0 \\
  S_1^*(\xi_{k_1,\ldots,k_n}), & i=1 \\
  S_i(\xi_{k_1,\ldots,k_n}), & 2\leq i\leq n \end{cases} \\
\sigma_-(e_i)\xi_{k_1,\ldots,k_n} & = & \begin{cases}
  -q^{(k_1+k_2+\cdots+k_n)/2}\xi_{k_1,\ldots,k_n}, & i=0 \\
  S_i(\xi_{k_1,\ldots,k_n}), & 1\leq i\leq n \end{cases}
\end{eqnarray*}
\end{prop}
\begin{proof}
Throughout this proof we denote by $a_j,b_j$, $j=0,\ldots,n$ the elements of
$C(B^{2n}_q)$ such that $e_j=(a_j,b_j)$ in formulae
(\ref{s2ngenerators5})--(\ref{s2ngenerators4}).

To calculate irreducible representations of $C(S^{2n}_{q,\beta})$ we
use exact sequence (\ref{s2nexact}). Thus $\irr(C(S^{2n}_{q,\beta}))$ is
the disjoint union of $\irr(J_{2n}\oplus J_{2n})$ and $\irr(C(\partial B^{2n}_q))$.

Since $J_{2n}\oplus0\cong\K$, this ideal has a unique irreducible representation,
whose extension to $C(S^{2n}_{q,\beta})$ we denote $\sigma_+$. If we identify
$C(S^{2n}_{q,\beta})$ with the algebra $C(B^{2n}_q)\oplus_\beta C(B^{2n}_q)$ of
the pull-back diagram (\ref{double}), then for $x\in J_{2n}$ we have
identification $e_j(x\oplus0)=(a_jx,0)$. Consequently $\sigma_+(e_j)=
\sigma(a_j)$, where $\sigma$ is the irreducible representation of $C(B^{2n}_q)$
from Proposition \ref{bevenirr}. Similarly, $\sigma_-(e_j)=\sigma(b_j)$.

Let $\rho$ be an irreducible representation of $C(\partial B^{2n}_q)$ and
let $\tilde{\rho}$ be its lift to $C(S^{2n}_{q,\beta})$. Then $\rho\circ\pi
\in\irr(C(B^{2n}_q))$ and $\rho\circ\pi(J_{2n})=\{0\}$. Thus there exist
$i$ and $\theta$ such that $\rho\circ\pi=\rho_i^\theta$, where $\rho_i^\theta$
is one of the representations of $C(B^{2n}_q)$ from Proposition
\ref{bevenirr}. But the quotient of $C(S^{2n}_{q,\beta})$ by the ideal
$J_{2n}\oplus0$ may be identified with $C(B^{2n}_q)$ in such a way that
the image  (under the natural surjection) of each $e_j$ is $a_j$. Thus
$\tilde{\rho}(e_j)=\rho_i^\theta(a_j)$ and a straightforward calculation
yields $\tilde{\rho}=\varrho_i^\theta$, as above.
\end{proof}

\end{document}